\documentclass[11pt]{article}
\usepackage{amsmath}
\usepackage{amssymb}
\usepackage{amsfonts}
\usepackage{eucal}
\usepackage{graphicx,url}
\usepackage{amssymb,amsmath, amsthm}
\usepackage{epstopdf}
\newtheorem{Theorem}{Theorem}[section]
\newtheorem{Proposition}[Theorem]{Proposition}

\newtheorem{Lemma}{Lemma}[section]
\newtheorem{Remark}{Remark}[section]
\pdfoutput=1
\begin{document}
\title{Large deviations for a damped telegraph process}
\author{Alessandro De Gregorio\thanks{Dipartimento di Scienze
Statistiche, Sapienza Universit\`{a} di Roma, Piazzale Aldo Moro
5, I-00185 Rome, Italy. e-mail:
\texttt{alessandro.degregorio@uniroma1.it}}\and Claudio
Macci\thanks{Dipartimento di Matematica, Universit\`a di Roma Tor
Vergata, Via della Ricerca Scientifica, I-00133 Rome, Italy.
e-mail: \texttt{macci@mat.uniroma2.it}}}
\date{}
\maketitle
\begin{abstract}
\noindent In this paper we consider a slight generalization of the
damped telegraph process in Di Crescenzo and Martinucci (2010). We
prove a large deviation principle for this process and an
asymptotic result for its level crossing probabilities (as the
level goes to infinity). Finally we compare our results with the
analogous well-known results for the standard telegraph process.
\end{abstract}

\section{Introduction}\label{sec:introduction}
The theory of large deviations gives an asymptotic computation of
small probabilities on exponential scale. Estimates based on large
deviations play a crucial role in resolving a variety of problems
in several fields. A part of these problems has interest in risk
theory and are solved by considering large deviation estimates for
some level crossing probabilities as, for instance, the ruin
probabilities for some insurance models or the overflow
probabilities for some queueing models.

In this paper we consider the damped telegraph process in Di
Crescenzo and Martinucci (2010) which is derived from the standard
telegraph process in Beghin et al. (2001); actually we have in
mind the case with drift (see Orsingher (1990) for the case
without drift). More precisely we consider the process
$\{D(t):t\geq 0\}$ which is a slight generalization of the one in
Di Crescenzo and Martinucci (2010) because an arbitrary
distribution for the random initial velocity is allowed (see eq.
\eqref{eq:arbitrary-distribution-for-initial-velocity} below). A
recent paper on large deviations for some telegraph processes is
De Gregorio and Macci (2012).

There is a wide literature on several versions of the telegraph
process, with applications; here we recall Mazza and Rulli\`{e}re
(2004) which illustrated an interesting link between the standard
telegraph process and the standard risk process in insurance (we
mean the \emph{compound Poisson model} in Section 5.3 in Rolski et
al. (1999), or the \emph{Cram\'{e}r-Lundberg model} in Section 1.1
in Embrechts et al. (1997)) with exponentially distributed claim
sizes. The results in this paper have interest for the asymptotic
behavior of some item modeled on a semi-Markov process (a wide
source of models can be found in Janssen and Manca (2006, 2007));
actually the random evolution of $\{D(t):t\geq 0\}$ is driven by a
continuous time Markov chain with two states and linearly
increasing switching rates, and therefore it is driven by a
particular non-homogeneous semi-Markov process.

A result in this paper concerns the probability that the process
$\{D(t):t\geq 0\}$ crosses the level $q$ on the infinite time
horizon $[0,\infty)$, i.e.
\begin{equation}\label{eq:LCP}
P(Q_D>q),\ \mathrm{where}\ Q_D:=\sup\{D(t):t\geq 0\}.
\end{equation}
Then, under a stability condition (see eq.
\eqref{eq:stability-condition} below), we prove that
\begin{equation}\label{eq:Lundberg-estimate-intro}
\lim_{q\to\infty}\frac{1}{q}\log P(Q_D>q)=-w_D
\end{equation}
for some $w_D>0$. The limit \eqref{eq:Lundberg-estimate-intro} is
proved by combining the \emph{large deviation principle} of
$\left\{\frac{D(t)}{t}:t\geq 0\right\}$ (as $t\to\infty$) proved
in this paper, and a quite general result in Duffy et al. (2003);
actually $w_D$ can be expressed in terms of the large deviation
rate function $I_D$ for $\left\{\frac{D(t)}{t}:t\geq 0\right\}$
(see eq. \eqref{eq:variational-formula} below).

We remark that the limit \eqref{eq:Lundberg-estimate-intro} has an
analogy with several results in the literature: here we recall
Duffy et al. (2003) cited above, Djehiche (1993) which provides a
result for risk processes with reserve dependent premium rate,
Lehtonen and Nyrhinen (1992a, 1992b) where the limit
\eqref{eq:Lundberg-estimate-intro} plays a crucial role for the
use of importance sampling technique in an estimation problem by
Monte Carlo simulations. In several cases the limit
\eqref{eq:Lundberg-estimate-intro} has a strict relationship with
some sharp exponential upper bounds for level crossing
probabilities, as for instance the well-known \emph{Lundberg
inequality} for random walks or L\'{e}vy processes (see e.g.
Theorem 5.1 in Asmussen (2003)). We also recall that in some cases
the Lundberg inequality can be seen as entropy estimate with an
interesting structure familiar from thermodynamics (see the
discussion in Martin-L\"{o}f (1986)). The only sharp upper bound
recalled in this paper concerns the standard telegraph process
(see Remark \ref{rem:sharpUBstandard} below).

We conclude with the outline of the paper. We start with some
preliminaries in Section \ref{sec:preliminaries}. In Section
\ref{sec:model} we present the damped telegraph process in this
paper. The results are presented in Section \ref{sec:results}.
Finally in Section \ref{sec:conclusions} we compare the results
obtained in this paper with the analogous well-known results for
the standard telegraph process, and we illustrate some open
problems.

\section{Preliminaries}\label{sec:preliminaries}
We start by recalling some basic definitions (see Dembo and
Zeitouni (1998), pages 4-5). Given a topological space
$\mathcal{Z}$ (here we always consider $\mathcal{Z}=\mathbb{R}$),
we say that a family of $\mathcal{Z}$-valued random variables
$\{Z(t):t>0\}$ satisfies the large deviation principle (LDP from
now on) with rate function $I$ if: the function $I:\mathcal{Z}\to
[0,\infty]$ is lower semi-continuous; the upper bound
$$\limsup_{t\to\infty}\frac{1}{t}\log P(Z(t)\in C)\leq-\inf_{x\in C}I(x)$$
holds for all closed sets $C$; the lower bound
$$\liminf_{t\to\infty}\frac{1}{t}\log P(Z(t)\in G)\geq-\inf_{x\in G}I(x)$$
holds for all open sets $G$. Moreover a rate function is said to
be good if all its level sets
$\{\{x\in\mathcal{Z}:I(x)\leq\eta\}:\eta\geq 0\}$ are compact.

Finally we recall Theorem 2.2 in Duffy et al. (2003). Here, for
simplicity, we present a slightly weaker version of the result;
more precisely, if we refer to the items in Duffy et al. (2003),
the functions $v$ and $a$ are defined by $v(t)=a(t)=t$ for all
$t>0$, and therefore we have $V=A=1$ and $h(t)=t$ for all $t>0$.

\begin{Proposition}\label{prop:DLS}
Assume that $\left\{\frac{X(t)}{t}:t>0\right\}$ satisfies the LDP
on $\mathbb{R}$ with rate function $I_X$ such that:\\
(i) $\inf_{x\geq 0}I_X(x)>0$;\\
(ii) there exists $y>0$ such that $\inf_{x\geq y}I_X(x)<\infty$;\\
(iii) the function $(0,\infty)\ni y\mapsto\inf_{x\geq y}I_X(x)$ is
continuous on the interior of the set upon which it is finite;\\
(iv) there exist $F>1$ and $K>0$ such that $\frac{1}{t}\log
P(X(t)>xt)\leq -x^F$ for all $t>0$ and for all $x>F$.\\
Then, if we set $Q_X^*:=\sup\{X(t):t\in\mathbb{N}\cup\{0\}\}$, we
have
$$\lim_{q\to\infty}\frac{1}{q}\log P(Q_X^*>q)=-w_X,\ where\ w_X:=\inf\left\{xI_X(1/x):x>0\right\}.$$
\end{Proposition}

\section{The damped telegraph process}\label{sec:model}
In this section we present the damped telegraph process studied in
this paper. We remark that it is a slight generalization of the
one studied by Di Crescenzo and Martinucci (2010); actually we
recover that model by setting $\alpha=\frac{1}{2}$.

We consider a random motion $\{D(t):t\geq 0\}$ on the real line
which starts at the origin and moves with a two-valued integrated
telegraph signal, i.e., for some $\lambda_1,\lambda_2,c_1,c_2>0$,
we have a rightward velocity $c_1$, a leftward velocity $-c_2$,
and the rates of the occurrences of velocity switches increase
linearly, i.e. they are $\lambda_1k$ and $\lambda_2k$ (for all
$k\geq 1$), respectively. More precisely we have
$$D(t):=\int_0^tV(s)ds,$$
where the velocity process $\{V(t):t\geq 0\}$ is defined by
$$V(t):=V(0)\left\{\frac{c_1-c_2}{2}+\frac{c_1+c_2}{2}\left\{1_{\{V(0)=c_1\}}-1_{\{V(0)=-c_2\}}\right\}(-1)^{N(t)}\right\},$$
and the random variable $V(0)$ is such that
$P(V(0)\in\{-c_2,c_1\})=1$. Moreover for the process $\{N(t):t\geq
0\}$ (which counts the number of changes of direction of
$\{D(t):t\geq 0\}$) we have $N(t):=\sum_{n\geq
1}1_{\{\tau_1+\cdots+\tau_n\leq t\}}$, where the random time
lengths $\{\tau_n:n\geq 1\}$ are conditionally independent given
$V(0)$, and the conditional distributions are the following:
$$\left.\begin{array}{ll}
\mathrm{if}\ V(0)=c_1,\ \mathrm{then}\ \left\{\begin{array}{ll}
\tau_{2k-1}\ \mbox{is exponentially distributed with mean}\ \frac{1}{\lambda_1k}\ (k\geq 1)\\
\tau_{2k}\ \mbox{is exponentially distributed with mean}\
\frac{1}{\lambda_2k}\ (k\geq 1);
\end{array}\right.\\
\mathrm{if}\ V(0)=-c_2,\ \mathrm{then}\ \left\{\begin{array}{ll}
\tau_{2k-1}\ \mbox{is exponentially distributed with mean}\ \frac{1}{\lambda_2k}\ (k\geq 1)\\
\tau_{2k}\ \mbox{is exponentially distributed with mean}\
\frac{1}{\lambda_1k}\ (k\geq 1).
\end{array}\right.
\end{array}\right.$$
Here we allow a general initial distribution of $V(0)$, i.e. we
set
\begin{equation}\label{eq:arbitrary-distribution-for-initial-velocity}
(P(V(0)=c_1),P(V(0)=-c_2))=(\alpha,1-\alpha)\ \mbox{for some}\
\alpha\in [0,1];
\end{equation}
as we shall see the results in this paper do not depend on the
value $\alpha$.

\begin{Remark}\label{rem:standard}
The process $\{D(t):t\geq 0\}$ is a suitable change of the
standard telegraph process $\{S(t):t\geq 0\}$ where the rates of
the occurrences of velocity switches are (the constant values)
$\lambda_1$ and $\lambda_2$, respectively.
\end{Remark}

We conclude with the generalization of Theorem 3.2 in Di Crescenzo
and Martinucci (2010), which can be recovered by setting
$\alpha=\frac{1}{2}$.

\begin{Lemma}\label{lem:Theorem-3.2-DM-generalization}
Let $t>0$ be arbitrarily fixed. Then we have
$$P(D(t)\in A)=\alpha e^{-\lambda_1t}1_A(c_1t)+(1-\alpha)e^{-\lambda_2t}1_A(-c_2t)+\int_Ap(x,t)1_{(-c_2t,c_1t)}(x)dx$$
for any Borel subset $A$ of $\mathbb{R}$, where
$$p(x,t)=\frac{e^{-\lambda_1\tau_*}e^{-\lambda_2(t-\tau_*)}[\lambda_1+\lambda_2
-\alpha\lambda_2e^{-\lambda_1\tau_*}-(1-\alpha)\lambda_1e^{-\lambda_2(t-\tau_*)}]}
{(c_1+c_2)[e^{-\lambda_2(t-\tau_*)}+e^{-\lambda_1\tau_*}(1-e^{-\lambda_2(t-\tau_*)})]^2}$$
and $\tau_*=\tau_*(x,t):=\frac{c_2t+x}{c_1+c_2}$.
\end{Lemma}
\noindent\emph{Proof.} It is immediate to check that we have the
probability masses equal to $\alpha e^{-\lambda_1t}$ and
$(1-\alpha)e^{-\lambda_2t}$ concentrated at the points $c_1t$ and
$-c_2t$, respectively. For the density on $(-c_2t,c_1t)$, we
follow the same lines of the procedures in Di Crescenzo and
Martinucci (2010) with some changes of notation (as far as that
reference is concerned, see eqs. (2.3) and (2.4), the successive
formulas in Section 3 and Remark 2.1, and the proof of Theorem
3.2): more precisely we have
$$p(x,t)=\alpha\{f(x,t|c_1)+b(x,t|c_1)\}+(1-\alpha)\{f(x,t|-c_2)+b(x,t|-c_2)\},$$
where
$$f(x,t|y):=\frac{\partial}{\partial x}P(D(t)\leq x,V(t)=c_1|V(0)=y)$$
and
$$b(x,t|y):=\frac{\partial}{\partial x}P(D(t)\leq x,V(t)=-c_2|V(0)=y)$$
(for $x\in(-c_2t,c_1t)$ and $y\in\{-c_2,c_1\}$). $\Box$

\section{Results}\label{sec:results}
Our aim is to prove large deviation results for the process
$\{D(t):t\geq 0\}$ presented in Section \ref{sec:model}. More
precisely we mean:
\begin{itemize}
\item \emph{Proposition \ref{prop:LDP}}: the LDP of
$\left\{\frac{D(t)}{t}:t>0\right\}$;
\item \emph{Proposition \ref{prop:DLS-consequence}}: an
asymptotic result (as $q\to\infty$) for the level crossing
probability in \eqref{eq:LCP}, i.e. the limit
\eqref{eq:Lundberg-estimate-intro} for some $w_D>0$, under the
stability condition
\begin{equation}\label{eq:stability-condition}
\lambda_2c_1-\lambda_1c_2<0
\end{equation}
which ensures that $D(t)$ goes to $-\infty$ as $t\to\infty$.
\end{itemize}

We remark that \eqref{eq:stability-condition} also ensures that
$S(t)$ in Remark \ref{rem:standard} goes to $-\infty$ as
$t\to\infty$; moreover \eqref{eq:stability-condition} is
equivalent to $\frac{c_1}{\lambda_1}<\frac{c_2}{\lambda_2}$, where
$\frac{1}{\lambda_1}$ is the mean of the random time intervals
where $S(t)$ moves with rightward velocity $c_1$, and
$\frac{1}{\lambda_2}$ is the mean of the random time intervals
where $S(t)$ moves with leftward velocity $-c_2$.

We start with the LDP of $\left\{\frac{D(t)}{t}:t>0\right\}$. The
proof is based on the same method used in the paper of Duffy and
Sapozhnikov (2008); actually, in that reference, Theorem 2 is
proved as a consequence of Theorems 3 and 4 which correspond to
\eqref{eq:forLB} and \eqref{eq:forUB} in this paper (see below),
respectively. We can consider this method because the random
variables $\left\{\frac{D(t)}{t}:t>0\right\}$ take values on a
compact set $[-c_2,c_1]$, and the LDP follows from, for example,
Theorem 4.1.11 in Dembo and Zeitouni (1998).

\begin{Proposition}\label{prop:LDP}
The family of random variables $\left\{\frac{D(t)}{t}:t>0\right\}$
satisfies the LDP with good rate function $I_D$ defined by
\begin{align*}
I_D(x):=&\left\{\begin{array}{ll}
\frac{|(\lambda_1+\lambda_2)x-(\lambda_2c_1-\lambda_1c_2)|}{c_1+c_2}&\
if\
x\in[-c_2,c_1]\\
\infty&\ otherwise
\end{array}\right.\\
=&\left\{\begin{array}{ll}
\frac{-(\lambda_1+\lambda_2)x+\lambda_2c_1-\lambda_1c_2}{c_1+c_2}&\
if\
x\in\left[-c_2,\frac{\lambda_2c_1-\lambda_1c_2}{\lambda_1+\lambda_2}\right]\\
\frac{(\lambda_1+\lambda_2)x-(\lambda_2c_1-\lambda_1c_2)}{c_1+c_2}&\
if\
x\in\left(\frac{\lambda_2c_1-\lambda_1c_2}{\lambda_1+\lambda_2},c_1\right]\\
\infty&\ otherwise.
\end{array}\right.
\end{align*}
\end{Proposition}
\noindent\emph{Proof.} We start with the following equality which
is often used throughout this proof:
\begin{equation}\label{eq:auxiliary-equality}
I_D(x)=\lambda_1\frac{c_2+x}{c_1+c_2}+\lambda_2\frac{c_1-x}{c_1+c_2}
+2\max\left\{-\lambda_2\frac{c_1-x}{c_1+c_2},-\lambda_1\frac{c_2+x}{c_1+c_2}\right\}\
\mbox{for all}\ x\in[-c_2,c_1].
\end{equation}
The equality \eqref{eq:auxiliary-equality} can be checked by
inspection. It is useful to distinguish the following three cases:
\begin{enumerate}
\item
$\lambda_2\frac{c_1-x}{c_1+c_2}<\lambda_1\frac{c_2+x}{c_1+c_2}$,
which is equivalent to
$x>\frac{\lambda_2c_1-\lambda_1c_2}{\lambda_1+\lambda_2}$;
\item
$\lambda_2\frac{c_1-x}{c_1+c_2}>\lambda_1\frac{c_2+x}{c_1+c_2}$,
which is equivalent to
$x<\frac{\lambda_2c_1-\lambda_1c_2}{\lambda_1+\lambda_2}$;
\item
$\lambda_2\frac{c_1-x}{c_1+c_2}=\lambda_1\frac{c_2+x}{c_1+c_2}$,
which is equivalent to
$x=\frac{\lambda_2c_1-\lambda_1c_2}{\lambda_1+\lambda_2}$.
\end{enumerate}
We also remark that we have $I_D(x)>0$ in the cases 1 and 2, and
$I_D(x)=0$ in the case 3.

We prove the LDP by checking the following asymptotic estimates
for all $x\in\mathbb{R}$:
\begin{equation}\label{eq:forLB}
\lim_{\varepsilon\to 0}\liminf_{t\to\infty}\frac{1}{t}\log
P\left(\frac{D(t)}{t}\in(x-\varepsilon,x+\varepsilon)\right)\geq-I_D(x);
\end{equation}
\begin{equation}\label{eq:forUB}
\lim_{\varepsilon\to 0}\limsup_{t\to\infty}\frac{1}{t}\log
P\left(\frac{D(t)}{t}\in(x-\varepsilon,x+\varepsilon)\right)\leq-I_D(x).
\end{equation}
We have the following cases.

\noindent $\bullet$ \emph{Case} $x\notin [-c_2,c_1]$. In this case
the proof of \eqref{eq:forLB} and \eqref{eq:forUB} is immediate
because we have $I_D(x)=\infty$. Actually \eqref{eq:forLB}
trivially holds; moreover, if we take $\varepsilon>0$ small enough
to have $x+\varepsilon<-c_2$ or $x-\varepsilon>c_1$, we have
$P\left(\frac{D(t)}{t}\in(x-\varepsilon,x+\varepsilon)\right)=0$
for all $t>0$, which yields \eqref{eq:forUB}.

\noindent $\bullet$ \emph{Case} $x\in (-c_2,c_1)$. Without loss of
generality we can take $\varepsilon>0$ small enough to have
$x-\varepsilon,x+\varepsilon\in(-c_2,c_1)$. Then there exists
$\tilde{z}=\tilde{z}(\varepsilon,t,x)\in(x-\varepsilon,x+\varepsilon)$
such that
\begin{equation}\label{eq:mean-value-theorem}
P\left(\frac{D(t)}{t}\in(x-\varepsilon,x+\varepsilon)\right)=\int_{(x-\varepsilon)t}^{(x+\varepsilon)t}p(y,t)dy=
\int_{(x-\varepsilon)}^{(x+\varepsilon)}p(zt,t)tdz=p(\tilde{z}t,t)t2\varepsilon.
\end{equation}
Moreover we remark that
$$\tau_*(\tilde{z}t,t)=\frac{c_2+\tilde{z}}{c_1+c_2}t\in\left(\frac{c_2+x-\varepsilon}{c_1+c_2}t,\frac{c_2+x+\varepsilon}{c_1+c_2}t\right)$$
and
$$t-\tau_*(\tilde{z}t,t)=t-\frac{c_2+\tilde{z}}{c_1+c_2}t=\frac{c_1-\tilde{z}}{c_1+c_2}t
\in\left(\frac{c_1-x-\varepsilon}{c_1+c_2}t,\frac{c_1-x+\varepsilon}{c_1+c_2}t\right).$$
Thus, by Lemma \ref{lem:Theorem-3.2-DM-generalization}, we have
\begin{equation}\label{eq:forLB-1}
p(\tilde{z}t,t)\geq\frac{e^{-\lambda_1\frac{c_2+x+\varepsilon}{c_1+c_2}t}e^{-\lambda_2\frac{c_1-x+\varepsilon}{c_1+c_2}t}
\left[\lambda_1+\lambda_2-\alpha\lambda_2e^{-\lambda_1\frac{c_2+x-\varepsilon}{c_1+c_2}t}-
(1-\alpha)\lambda_1e^{-\lambda_2\frac{c_1-x-\varepsilon}{c_1+c_2}t}\right]}
{(c_1+c_2)\left[e^{-\lambda_2\frac{c_1-x-\varepsilon}{c_1+c_2}t}+e^{-\lambda_1\frac{c_2+x-\varepsilon}{c_1+c_2}t}\right]^2}
\end{equation}
and
\begin{equation}\label{eq:forUB-1}
p(\tilde{z}t,t)\leq\frac{e^{-\lambda_1\frac{c_2+x-\varepsilon}{c_1+c_2}t}e^{-\lambda_2\frac{c_1-x-\varepsilon}{c_1+c_2}t}
\left[\lambda_1+\lambda_2-\alpha\lambda_2e^{-\lambda_1\frac{c_2+x+\varepsilon}{c_1+c_2}t}-
(1-\alpha)\lambda_1e^{-\lambda_2\frac{c_1-x+\varepsilon}{c_1+c_2}t}\right]}
{(c_1+c_2)\left[e^{-\lambda_2\frac{c_1-x+\varepsilon}{c_1+c_2}t}+e^{-\lambda_1\frac{c_2+x+\varepsilon}{c_1+c_2}t}
\left(1-e^{-\lambda_2\frac{c_1-x-\varepsilon}{c_1+c_2}t}\right)\right]^2}.
\end{equation}

\noindent\emph{Proof of \eqref{eq:forLB} for $x\in(-c_2,c_1)$}.
Firstly, by \eqref{eq:mean-value-theorem} and \eqref{eq:forLB-1},
we have
\begin{align*}
\liminf_{t\to\infty}\frac{1}{t}\log
P\left(\frac{D(t)}{t}\in(x-\varepsilon,x+\varepsilon)\right)=&\liminf_{t\to\infty}\frac{1}{t}\log p(\tilde{z}t,t)\\
\geq&-\lambda_1\frac{c_2+x+\varepsilon}{c_1+c_2}-\lambda_2\frac{c_1-x+\varepsilon}{c_1+c_2}\\
&-2\limsup_{t\to\infty}\frac{1}{t}\log\left(e^{-\lambda_2\frac{c_1-x-\varepsilon}{c_1+c_2}t}
+e^{-\lambda_1\frac{c_2+x-\varepsilon}{c_1+c_2}t}\right);
\end{align*}
then, by considering Lemma 1.2.15 in Dembo and Zeitouni (1998) for
the last term, we obtain
\begin{align*}
\liminf_{t\to\infty}\frac{1}{t}\log
P\left(\frac{D(t)}{t}\in(x-\varepsilon,x+\varepsilon)\right)
\geq&-\lambda_1\frac{c_2+x+\varepsilon}{c_1+c_2}-\lambda_2\frac{c_1-x+\varepsilon}{c_1+c_2}\\
&-2\max\left\{-\lambda_2\frac{c_1-x-\varepsilon}{c_1+c_2},-\lambda_1\frac{c_2+x-\varepsilon}{c_1+c_2}\right\};
\end{align*}
finally we get \eqref{eq:forLB} by letting $\varepsilon$ go to
zero and by taking into account \eqref{eq:auxiliary-equality}.\\
\noindent\emph{Proof of \eqref{eq:forUB} for $x\in(-c_2,c_1)$}. We
introduce the symbol $A(t,\varepsilon)$ for the denominator in the
right hand side of \eqref{eq:forUB-1}:
$$A(t,\varepsilon):=(c_1+c_2)\left[e^{-\lambda_2\frac{c_1-x+\varepsilon}{c_1+c_2}t}+e^{-\lambda_1\frac{c_2+x+\varepsilon}{c_1+c_2}t}
\left(1-e^{-\lambda_2\frac{c_1-x-\varepsilon}{c_1+c_2}t}\right)\right]^2.$$
Then, by \eqref{eq:mean-value-theorem} and \eqref{eq:forUB-1}, we
have
\begin{align*}
\limsup_{t\to\infty}\frac{1}{t}\log
P\left(\frac{D(t)}{t}\in(x-\varepsilon,x+\varepsilon)\right)=&\limsup_{t\to\infty}\frac{1}{t}\log p(\tilde{z}t,t)\\
\leq&-\lambda_1\frac{c_2+x-\varepsilon}{c_1+c_2}-\lambda_2\frac{c_1-x-\varepsilon}{c_1+c_2}-\liminf_{t\to\infty}\frac{1}{t}\log
A(t,\varepsilon);
\end{align*}
moreover, if we take into account that
$$A(t,\varepsilon)\geq
\left\{\begin{array}{ll}
(c_1+c_2)e^{-2\lambda_2\frac{c_1-x+\varepsilon}{c_1+c_2}t}&
\ \mathrm{if}\ -\lambda_2\frac{c_1-x+\varepsilon}{c_1+c_2}\geq-\lambda_1\frac{c_2+x+\varepsilon}{c_1+c_2}\\
(c_1+c_2)e^{-2\lambda_1\frac{c_2+x+\varepsilon}{c_1+c_2}t}\left(1-e^{-\lambda_2\frac{c_1-x-\varepsilon}{c_1+c_2}t}\right)^2&\
\mathrm{if}\
-\lambda_2\frac{c_1-x+\varepsilon}{c_1+c_2}<-\lambda_1\frac{c_2+x+\varepsilon}{c_1+c_2},
\end{array}\right.$$
we get
$$\liminf_{t\to\infty}\frac{1}{t}\log A(t,\varepsilon)\geq
\max\left\{-2\lambda_2\frac{c_1-x+\varepsilon}{c_1+c_2},-2\lambda_1\frac{c_2+x+\varepsilon}{c_1+c_2}\right\}.$$
Then we obtain
\begin{align*}
\limsup_{t\to\infty}\frac{1}{t}\log
P\left(\frac{D(t)}{t}\in(x-\varepsilon,x+\varepsilon)\right)
\leq&-\lambda_1\frac{c_2+x-\varepsilon}{c_1+c_2}-\lambda_2\frac{c_1-x-\varepsilon}{c_1+c_2}\\
&-2\max\left\{-\lambda_2\frac{c_1-x+\varepsilon}{c_1+c_2},-\lambda_1\frac{c_2+x+\varepsilon}{c_1+c_2}\right\};
\end{align*}
finally we get \eqref{eq:forUB} by letting $\varepsilon$ go to
zero and by taking into account \eqref{eq:auxiliary-equality}.

\noindent $\bullet$ \emph{Case} $x=c_1$. It is similar to the case
$x\in(-c_2,c_1)$ with suitable changes; roughly speaking we often
have to consider $x=c_1$ in place of
$x+\varepsilon=c_1+\varepsilon$. We start with the analogous of
\eqref{eq:mean-value-theorem}, \eqref{eq:forLB-1} and
\eqref{eq:forUB-1}:
$$P\left(\frac{D(t)}{t}\in(c_1-\varepsilon,c_1+\varepsilon)\right)=p(\tilde{z}t,t)t\varepsilon$$
for some
$\tilde{z}=\tilde{z}(\varepsilon,t,c_1)\in(c_1-\varepsilon,c_1)$;
moreover, since
$$\tau_*(\tilde{z}t,t)\in\left(\frac{c_2+c_1-\varepsilon}{c_1+c_2}t,t\right)\
\mbox{and}\
t-\tau_*(\tilde{z}t,t)\in\left(0,\frac{\varepsilon}{c_1+c_2}t\right),$$
we have
$$p(\tilde{z}t,t)\geq\frac{e^{-\lambda_1t}e^{-\lambda_2\frac{\varepsilon}{c_1+c_2}t}
\left[\lambda_1+\lambda_2-\alpha\lambda_2e^{-\lambda_1\frac{c_2+c_1-\varepsilon}{c_1+c_2}t}-(1-\alpha)\lambda_1\right]}
{(c_1+c_2)\left[1+e^{-\lambda_1\frac{c_2+c_1-\varepsilon}{c_1+c_2}t}\right]^2}$$
and
$$p(\tilde{z}t,t)\leq\frac{e^{-\lambda_1\frac{c_2+c_1-\varepsilon}{c_1+c_2}t}
\left[\lambda_1+\lambda_2-\alpha\lambda_2e^{-\lambda_1t}-(1-\alpha)\lambda_1e^{-\lambda_2\frac{\varepsilon}{c_1+c_2}t}\right]}
{(c_1+c_2)\left[e^{-\lambda_2\frac{\varepsilon}{c_1+c_2}t}\right]^2}.$$
Thus
$$\lim_{\varepsilon\to 0}\liminf_{t\to\infty}\frac{1}{t}\log
P\left(\frac{D(t)}{t}\in(c_1-\varepsilon,c_1+\varepsilon)\right)\geq
\lim_{\varepsilon\to
0}-\lambda_1-\lambda_2\frac{\varepsilon}{c_1+c_2}=-\lambda_1=-I_D(c_1)$$
and
\begin{align*}
\lim_{\varepsilon\to 0}\limsup_{t\to\infty}\frac{1}{t}\log
P\left(\frac{D(t)}{t}\in(c_1-\varepsilon,c_1+\varepsilon)\right)
&\leq\lim_{\varepsilon\to
0}-\lambda_1\frac{c_2+c_1-\varepsilon}{c_1+c_2}+2\lambda_2\frac{\varepsilon}{c_1+c_2}\\
&=-\lambda_1=-I_D(c_1).
\end{align*}

\noindent $\bullet$ \emph{Case} $x=-c_2$. We argue as for the case
$x=c_1$. Thus we proceed as for the case $x\in(-c_2,c_1)$ with
suitable changes; roughly speaking we often have to consider
$x=-c_2$ in place of $x-\varepsilon=-c_2-\varepsilon$. We have
$$P\left(\frac{D(t)}{t}\in(-c_2-\varepsilon,-c_2+\varepsilon)\right)=p(\tilde{z}t,t)t\varepsilon$$
for some
$\tilde{z}=\tilde{z}(\varepsilon,t,-c_2)\in(-c_2,-c_2+\varepsilon)$;
moreover, since
$$\tau_*(\tilde{z}t,t)\in\left(0,\frac{\varepsilon}{c_1+c_2}t\right)\
\mbox{and}\
t-\tau_*(\tilde{z}t,t)\in\left(\frac{c_1+c_2-\varepsilon}{c_1+c_2}t,t\right),$$
we have
$$p(\tilde{z}t,t)\geq\frac{e^{-\lambda_1\frac{\varepsilon}{c_1+c_2}t}e^{-\lambda_2t}
\left[\lambda_1+\lambda_2-\alpha\lambda_2-(1-\alpha)\lambda_1e^{-\lambda_2\frac{c_1+c_2-\varepsilon}{c_1+c_2}t}\right]}
{(c_1+c_2)\left[e^{-\lambda_2\frac{c_1+c_2-\varepsilon}{c_1+c_2}t}+1\right]^2}$$
and
$$p(\tilde{z}t,t)\leq\frac{e^{-\lambda_2\frac{c_1+c_2-\varepsilon}{c_1+c_2}t}
\left[\lambda_1+\lambda_2-\alpha\lambda_2e^{-\lambda_1\frac{\varepsilon}{c_1+c_2}t}-(1-\alpha)\lambda_1e^{-\lambda_2t}\right]}
{(c_1+c_2)\left[e^{-\lambda_2t}+e^{-\lambda_1\frac{\varepsilon}{c_1+c_2}t}\left(1-e^{-\lambda_2\frac{c_1+c_2-\varepsilon}{c_1+c_2}t}\right)\right]^2}.$$
Thus
$$\lim_{\varepsilon\to 0}\liminf_{t\to\infty}\frac{1}{t}\log
P\left(\frac{D(t)}{t}\in(-c_2-\varepsilon,-c_2+\varepsilon)\right)\geq
\lim_{\varepsilon\to
0}-\lambda_1\frac{\varepsilon}{c_1+c_2}-\lambda_2=-\lambda_2=-I_D(-c_2)$$
and
\begin{align*}
\lim_{\varepsilon\to 0}\limsup_{t\to\infty}\frac{1}{t}\log
P\left(\frac{D(t)}{t}\in(c_1-\varepsilon,c_1+\varepsilon)\right)
&\leq\lim_{\varepsilon\to
0}-\lambda_2\frac{c_1+c_2-\varepsilon}{c_1+c_2}-2\max\left\{-\lambda_2,-\lambda_1\frac{\varepsilon}{c_1+c_2}\right\}\\
&=-\lambda_2=-I_D(-c_2).\ \Box\\
\end{align*}

Now we prove an asymptotic result (as $q\to\infty$) for the level
crossing probabilities in \eqref{eq:LCP}. This result will be
proved by applying Proposition \ref{prop:DLS} together with
Proposition \ref{prop:LDP}.

\begin{Proposition}\label{prop:DLS-consequence}
Assume that $\lambda_2c_1-\lambda_1c_2<0$. Then we have
$\lim_{q\to\infty}\frac{1}{q}\log
P(Q_D>q)=-\frac{\lambda_1}{c_1}$.
\end{Proposition}
\noindent\emph{Proof.} We want to apply Proposition
\ref{prop:DLS}; therefore we define
$Q_D^*=\sup\{D(t):t\in\mathbb{N}\cup\{0\}\}$. We remark that
$Q_D^*\leq Q_D\leq Q_D^*+c_1$ (obviously the two inequalities turn
into equalities if and only if $Q_D^*=\infty$; however we have
$P(Q_D^*=\infty)=0$ by the hypothesis
$\lambda_2c_1-\lambda_1c_2<0$); then it suffices to show that
\begin{equation}\label{eq:dt-restriction}
\lim_{q\to\infty}\frac{1}{q}\log P(Q_D^*>q)=-\frac{\lambda_1}{c_1}
\end{equation}
because, in such a case, we would have
$$\liminf_{q\to\infty}\frac{1}{q}\log P(Q_D>q)\geq\liminf_{q\to\infty}\frac{1}{q}\log P(Q_D^*>q)=-\frac{\lambda_1}{c_1}$$
and
$$\limsup_{q\to\infty}\frac{1}{q}\log P(Q_D>q)\leq\limsup_{q\to\infty}\frac{1}{q}\log P(Q_D^*>q-c_1)=-\frac{\lambda_1}{c_1},$$
which yield $\lim_{q\to\infty}\frac{1}{q}\log
P(Q_D>q)=-\frac{\lambda_1}{c_1}$.

The limit \eqref{eq:dt-restriction} can be proved by applying
Proposition \ref{prop:DLS} with $\{X(t):t\geq 0\}=\{D(t):t\geq
0\}$ and, by Proposition \ref{prop:LDP}, with $I_X=I_D$. We remark
that we have $\inf\{I_D(x):x\geq y\}=I_D(y)$ for all $y>0$ by the
hypothesis $\lambda_2c_1-\lambda_1c_2<0$. Thus hypotheses (i),
(ii) and (iii) in Proposition \ref{prop:DLS} trivially hold.
Moreover hypothesis (iv) in Proposition \ref{prop:DLS} holds with
$K=c_1$; actually, for all $t>0$, we have $P(D(t)>xt)=0$ for all
$x>c_1$ because $P(D(t)\in[-c_2t,c_1t])=1$. Then, by Proposition
\ref{prop:DLS}, we have $\lim_{q\to\infty}\frac{1}{q}\log
P(Q_D^*>q)=-w_D$, where
\begin{equation}\label{eq:variational-formula}
w_D:=\inf\left\{xI_D(1/x):x>0\right\}.
\end{equation}
We conclude with the computation of the infimum. We have
$$xI_D(1/x)=x\frac{\frac{\lambda_1+\lambda_2}{x}-(\lambda_2c_1-\lambda_1c_2)}{c_1+c_2}
=\frac{\lambda_1+\lambda_2-(\lambda_2c_1-\lambda_1c_2)x}{c_1+c_2}$$
for $x\geq 1/c_1$ (which is equivalent to $0\leq 1/x\leq c_1$);
then, again by the hypothesis $\lambda_2c_1-\lambda_1c_2<0$, the
infimum is attained at $x=1/c_1$, and we have
\begin{equation}\label{eq:Lundberg-parameter}
w_D=\frac{\lambda_1+\lambda_2-\frac{\lambda_2c_1-\lambda_1c_2}{c_1}}{c_1+c_2}
=\frac{(\lambda_1+\lambda_2)c_1-(\lambda_2c_1-\lambda_1c_2)}{c_1(c_1+c_2)}
=\frac{\lambda_1(c_1+c_2)}{c_1(c_1+c_2)}=\frac{\lambda_1}{c_1}.\
\Box
\end{equation}

\section{Conclusions}\label{sec:conclusions}
In this section we compare the results obtained for $\{D(t):t\geq
0\}$ with the well-known analogous results for $\{S(t):t\geq 0\}$
in Remark \ref{rem:standard}. In particular we show that, as one
expects because of the damping effect, the convergence at zero of
some rare events concerning $\{D(t):t\geq 0\}$ is faster than the
convergence of the analogous events concerning $\{S(t):t\geq 0\}$.
Finally we illustrate some open problems.

\subsection{Comparison between $\{D(t):t\geq 0\}$ and $\{S(t):t\geq 0\}$}
We start by recalling the analogous of Proposition \ref{prop:LDP}.
Here we refer to Macci (2009) but one could refer to Ney and
Nummelin (1987a, 1987b, 1987c) which concern the more general
setting of Markov additive processes.

\begin{Proposition}\label{prop:LDP-standard}
The family of random variables $\left\{\frac{S(t)}{t}:t>0\right\}$
satisfies the LDP with good rate function $I_S$ defined by
$$I_S(x):=\left\{\begin{array}{ll}
\left(\sqrt{\lambda_1\frac{x+c_2}{c_1+c_2}}-\sqrt{\lambda_2\frac{c_1-x}{c_1+c_2}}\right)^2&\ if\ x\in[-c_2,c_1]\\
\infty&\ otherwise.
\end{array}\right.$$
\end{Proposition}
\noindent\emph{Proof.} See Subsection 3.1 in Macci (2009); the
rate function $I_S$ coincides with $\kappa^*$ in that reference.
$\Box$\\

Now we recall the analogous of \ref{prop:DLS-consequence}, i.e.
the asymptotic result (as $q\to\infty$) for the level crossing
probabilities
$$P(Q_S>q),\ \mathrm{where}\ Q_S:=\sup\{S(t):t\geq 0\}.$$

\begin{Proposition}\label{prop:DLS-consequence-standard}
Assume that $\lambda_2c_1-\lambda_1c_2<0$. Then we have
$\lim_{q\to\infty}\frac{1}{q}\log
P(Q_S>q)=-\frac{\lambda_1c_2-\lambda_2c_1}{c_1c_2}$.
\end{Proposition}
\noindent\emph{Proof.} See Proposition 2.1 and Remark 2.2 in Macci
(2009), where $\frac{\lambda_1c_2-\lambda_2c_1}{c_1c_2}$ coincides
with $w$ in Subsection 3.1 in Macci (2009) (in the case where
$\mathbf{(H1)}$ holds). $\Box$

\begin{Remark}\label{rem:sharpUBstandard}
It is also known that we have a sharp upper bound for $P(Q_S>q)$;
more precisely (see e.g. Remark 2.3 in Macci (2009) which concerns
a more general Markov additive process) there exists $m>0$ such
that $P(Q_S>q)\leq
me^{-q\frac{\lambda_1c_2-\lambda_2c_1}{c_1c_2}}$ for all $q>0$.
\end{Remark}

\begin{figure}[ht]
\begin{center}
\includegraphics[angle=0,width=0.90\textwidth]{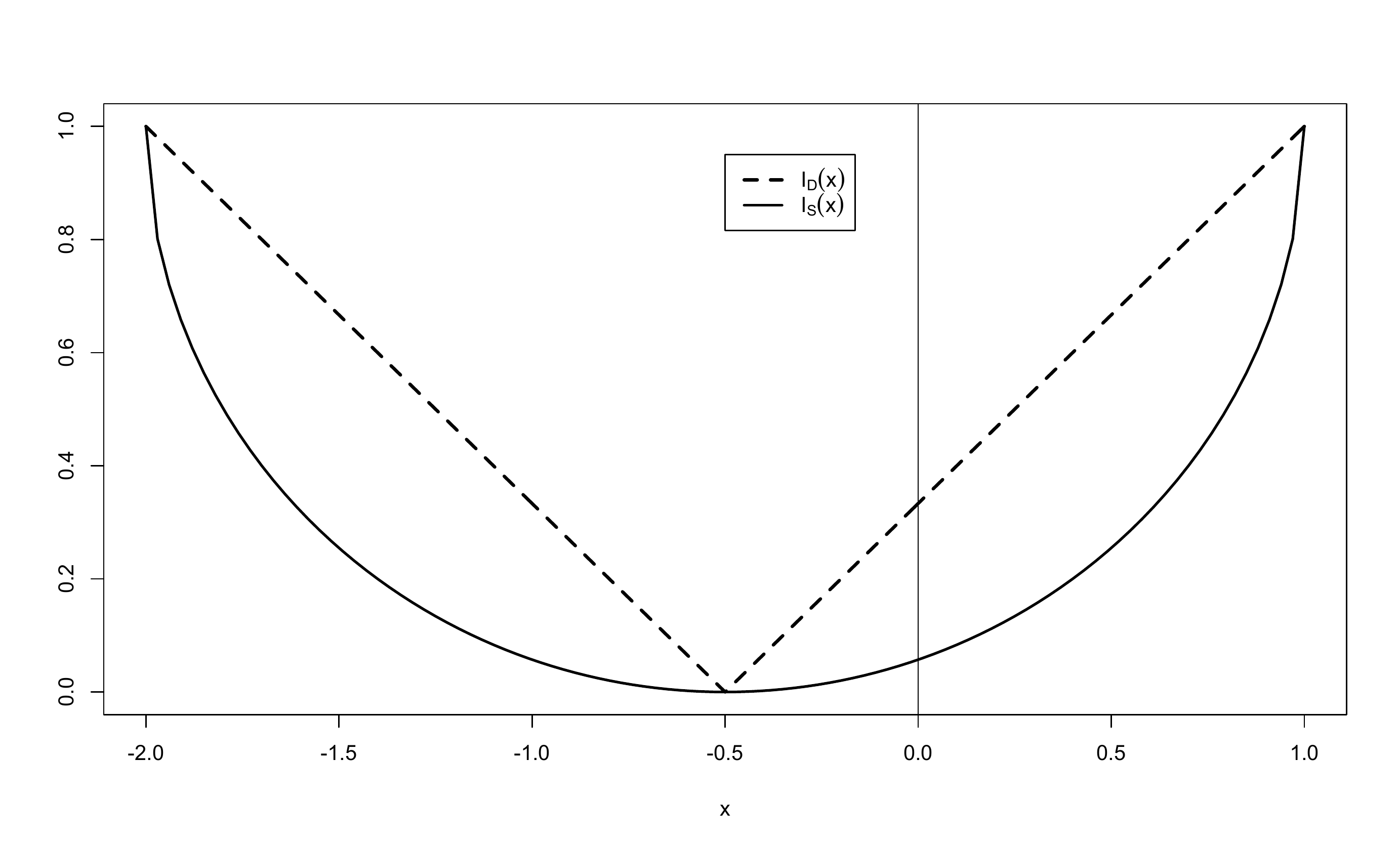}
\caption{The rate functions $I_D$ and $I_S$ in $[-c_2,c_1]$ for
$c_1=1$, $c_2=2$ and $\lambda_1=\lambda_2=1$.}\label{fig1}
\end{center}
\end{figure}

It is interesting to compare the rate function $I_D$ in
Proposition \ref{prop:LDP} and the rate function $I_S$ in
Proposition \ref{prop:LDP-standard}. Then we have the following
situation (see Figure \ref{fig1}).

\begin{itemize}
\item Both $I_D$ and $I_S$ uniquely vanish at
$x=\frac{\lambda_2c_1-\lambda_1c_2}{\lambda_1+\lambda_2}$;
therefore both $\frac{D(t)}{t}$ and $\frac{S(t)}{t}$ converge to
$\frac{\lambda_2c_1-\lambda_1c_2}{\lambda_1+\lambda_2}$ as
$t\to\infty$.
\item $I_D(x)=I_S(x)=\infty$ for all $x\notin [-c_2,c_1]$; actually
both $\left\{\frac{D(t)}{t}:t>0\right\}$ and
$\left\{\frac{S(t)}{t}:t>0\right\}$ are families of random
variables taking values on the closed set $[-c_2,c_1]$.
\item $I_D(-c_2)=I_S(-c_2)=\lambda_2$ and $I_D(c_1)=I_S(c_1)=\lambda_1$;
actually both cases $x=c_1$ and $x=-c_2$ concern the occurrence of
event \emph{no changes of direction} (and this event has the same
probability for both $\{D(t):t\geq 0\}$ and $\{S(t):t\geq 0\}$).
\item $I_D(x)>I_S(x)$ for all
$x\in(-c_2,c_1)\setminus\left\{\frac{\lambda_2c_1-\lambda_1c_2}{\lambda_1+\lambda_2}\right\}$.
\end{itemize}

As a consequence of the last statement we can say that, roughly
speaking, for any nonempty measurable set
$A\subset(-c_2,c_1)\setminus\left\{\frac{\lambda_2c_1-\lambda_1c_2}{\lambda_1+\lambda_2}\right\}$,
$P\left(\frac{D(t)}{t}\in A\right)$ converges to 0 faster than
$P\left(\frac{S(t)}{t}\in A\right)$ (as $t\to\infty$).

We remark that one can provide an alternative proof of Proposition
\ref{prop:DLS-consequence-standard} by following the same lines of
the proof of Proposition \ref{prop:DLS-consequence}. More
precisely one can check that $\lim_{q\to\infty}\frac{1}{q}\log
P(Q_S>q)=-w_S$, where
\begin{equation}\label{eq:variational-formula-standard}
w_S:=\inf\left\{xI_S(1/x):x>0\right\},
\end{equation}
and the equality
\begin{equation}\label{eq:Lundberg-parameter-standard}
w_S=\frac{\lambda_1c_2-\lambda_2c_1}{c_1c_2}.
\end{equation}
Thus the inequality $w_S\leq w_D$ is a straightforward consequence
of \eqref{eq:variational-formula},
\eqref{eq:variational-formula-standard} and the above detailed
inequality between the rate functions, i.e. $I_D(x)\geq I_S(x)$
for all $x\in\mathbb{R}$. However we can easily check the strict
inequality noting that
$$w_S=\frac{\lambda_1c_2-\lambda_2c_1}{c_1c_2}<\frac{\lambda_1c_2}{c_1c_2}=w_D$$
by \eqref{eq:Lundberg-parameter-standard} and
\eqref{eq:Lundberg-parameter}; thus, roughly speaking, $P(Q_D>q)$
converges to 0 faster than $P(Q_S>q)$ (as $q\to\infty$).

\subsection{Open problems}
The first open problem concerns moderate deviations. More
precisely, for some $\sigma^2\in(0,\infty)$, we should have the
following bounds for each $\{a_t:t>0\}$ such that $a_t\to 0$ and
$ta_t\to\infty$ (as $t\to\infty$):
$$\limsup_{t\to\infty}a_t\log P\left(\sqrt{\frac{a_t}{t}}(D(t)-\mathbb{E}[D(t)])\in C\right)\leq-\inf_{x\in C}\frac{x^2}{2\sigma^2}$$
for all closed sets $C$, and
$$\liminf_{t\to\infty} a_t\log P\left(\sqrt{\frac{a_t}{t}}(D(t)-\mathbb{E}[D(t)])\in G\right)\geq-\inf_{x\in G}\frac{x^2}{2\sigma^2}$$
for all open sets $G$. An analogue result for multivariate
centered random walks is Theorem 3.7.1 in Dembo and Zeitouni
(1998). We remark that, when one has the LDP with a convex a
regular rate function $I$ which uniquely vanishes at some point
$x_0$, the value $\sigma^2$ in the statement of the moderate
deviation result typically coincides with the inverse of
$I^{\prime\prime}(x_0)$. An interesting issue of this open problem
is that we cannot have this situation; in fact the rate function
$I$, i.e. $I_D$ in Proposition \ref{prop:LDP}, uniquely vanishes
at $x_0=\frac{\lambda_2c_1-\lambda_1c_2}{\lambda_1+\lambda_2}$ and
is not differentiable in $x_0$.

Another open problem concerns the case where the holding times are
heavy tailed distributed (and not exponentially distributed as
happens in the models studied in this paper). For instance one
could consider heavy tailed Weibull distributed holding times. In
this case, at least when all the holding times are equally
distributed, the LDP for the model without damping effect can be
obtained as a consequence of Theorem 2 in Duffy and Sapozhnikov
(2008). On the contrary we cannot say what happens for the damped
model.


\begin{thebibliography}{spc}
\bibitem{}
Asmussen, S. (2003). Applied Probability and Queues. Second
Edition. Springer-Verlag, New York.
\bibitem{}
Beghin, L., Nieddu, L., Orsingher, E. (2001). Probabilistic
analysis of the telegrapher's process with drift by means of
relativistic transformations. J. Appl. Math. Stochastic Anal. 14,
11--25.
\bibitem{}
De Gregorio, A., Macci, C. (2012). Large deviation principles for
telegraph processes. Statist. Probab. Lett. 82, 1874--1882.
\bibitem{}
Dembo, A., Zeitouni, O. (1998). Large Deviations Techniques and
Applications. Second Edition. Springer, New York.
\bibitem{}
Di Crescenzo, A., Martinucci, B. (2010). A damped telegraph random
process with logistic stationary distribution. J. Appl. Probab.
47, 84--96.
\bibitem{}
Djehiche, B. (1993). A large deviation estimate for ruin
probabilities. Scand. Actuar. J. 1993, no. 1, 42--59.
\bibitem{}
Duffy, K., Lewis, J.T., Sullivan, W.G. (2003). Logarithmic
asymptotics for the supremum of a stochastic processes. Ann. Appl.
Probab. 13, 430--445.
\bibitem{}
Duffy, K., Sapozhnikov, A. (2008). The large deviation principle
for the on-off Weibull sojourn process. J. Appl. Probab. 45,
107--117.
\bibitem{}
Embrechts, P., Kl\"{u}ppelberg, C., Mikosch, T. (1997). Modelling
Extremal Events. Springer, Berlin-Heidelberg.
\bibitem{}
Janssen, J., Manca, R. (2006). Applied Semi-Markov Processes. New
York, Springer.
\bibitem{}
Janssen, J., Manca, R. (2007). Semi-Markov Risk Models for
Finance, Insurance and Reliability. New York, Springer.
\bibitem{}
Lehtonen, T., Nyrhinen, H. (1992a). Simulating level crossing
probabilities by importance sampling. Adv. in Appl. Probab. 24,
858--874.
\bibitem{}
Lehtonen, T., Nyrhinen, H. (1992b). On asymptotically efficient
simulation of ruin probabilities in a Markovian environment.
Scand. Actuarial J. 1992, n. 1, 60--75.
\bibitem{}
Macci, C. (2009). Convergence of large deviation rates based on a
link between wave governed random motions and ruin processes.
Statist. Probab. Lett. 79, 255--263.
\bibitem{}
Martin-L\"{o}f, A. (1986). Entropy, a useful concept in risk
theory. Scand. Actuar. J. 1986, no. 3-4, 223--235.
\bibitem{}
Mazza, C., Rulli\`{e}re, D. (2004). A link between wave governed
random motions and ruin processes. Insurance Math. Econom. 35,
205--222.
\bibitem{}
Ney, P., Nummelin, E. (1987a). Markov additive processes I,
eigenvalue properties and limit theorems. Ann. Probab. 15,
561--592.
\bibitem{}
Ney, P., Nummelin, E. (1987b). Markov additive processes II, large
deviations. Ann. Probab. 15, 593--609.
\bibitem{}
Ney, P., Nummelin, E. (1987c). Markov additive processes: large
deviations for the continuous time case. In: Prohorov, Y.V.,
Statulevicius, V.A., Sazonov, V.V., Grigelionis, B. (Eds.),
Probability Theory and Mathematical Statistics (Vol. II). VNU Sci.
Press, Utrecht, pp. 377--389.
\bibitem{}
Orsingher, E. (1990). Probability law, flow function, maximum
distribution of wave governed random motions and their connections
with Kirchoff's laws. Stochastic Process. Appl. 34, 49--66.
\bibitem{}
Rolski, T., Schmidli, H., Schmidt, V., Teugels, J. (1999).
Stochastic Processes for Insurance and Finance. John Wiley and
Sons, Chichester.
\end{thebibliography}
\end{document}